\documentclass[onecolumn]{svjour3}
\usepackage{graphicx}
\usepackage{subfigure,graphicx}
\usepackage{amsmath}
\usepackage[]{placeins}
\usepackage{amsfonts}
\usepackage{float}
\usepackage{amssymb}
\usepackage{mathtools}
\usepackage{tensor}
\usepackage{enumitem}
\def\bfs{{\bf s}}
\def\bfu{{\bf u}}

\def\bfx{{\bf x}}

\def\bfE{{\bf E}}

\def\bfI{{\bf I}}

\def\bfN{{\bf N}}

\def\bfS{{\bf S}}

\def\bfX{{\bf X}}

\def\e0{\varepsilon_0}

\def\bfsig{\mbox{\boldmath $\sigma$}}

\newcommand\sts{\sigma_{\texttt{ts}}}
\newcommand\scs{\sigma_{\texttt{cs}}}
\newcommand\shs{\sigma_{\texttt{hs}}}
\newcommand\sbs{\sigma_{\texttt{bs}}}
\newcommand\sss{\sigma_{\texttt{ss}}}
\renewcommand\e{\varepsilon}

\newcommand\wc{\rightharpoonup}

\long\def\symbolfootnote[#1]#2{\begingroup%
\def\thefootnote{\fnsymbol{footnote}}\footnote[#1]{#2}\endgroup}

\begin{document}

%\journalname{}
\titlerunning{A variational formulation of Griffith phase-field fracture with material strength}

\title{A Variational Formulation of Griffith Phase-Field Fracture \\ with Material Strength}

\author{C. J. Larsen \and J. E. Dolbow \and  O. Lopez-Pamies}

\institute{
           Christopher J. Larsen \at Department of Mathematical Sciences, Worcester Polytechnic Institute, Worcester, MA 01609-2280, USA\\
           \email{cjlarsen@wpi.edu}\vspace{0.1cm}
           \and
           John E. Dolbow \at Department of Mechanical Engineering, Duke University, Durham, NC 27708, USA \\
           \email{john.dolbow@duke.edu}\vspace{0.1cm}
           \and
           Oscar Lopez-Pamies \at Department of Civil and Environmental Engineering, University of Illinois, Urbana--Champaign, IL 61801-2352, USA \\
           \email{pamies@illinois.edu}
           }

\date{}
\maketitle

\begin{abstract}

In this expository Note, it is shown that the Griffith phase-field theory of fracture accounting for material strength originally introduced by Kumar, Francfort, and Lopez-Pamies (J Mech Phys Solids 112, 523--551, 2018) in the form of PDEs can be recast as a variational theory. In particular, the solution pair $(\bfu,v)$ defined by the PDEs for the displacement field $\bfu$ and the phase field $v$  is shown to correspond to the fields that minimize separately two different functionals, much like the solution pair $(\bfu,v)$ defined by the original phase-field theory of fracture without material strength implemented in terms of alternating minimization. The merits of formulating a complete theory of fracture nucleation and propagation via such a variational approach --- in terms of the minimization of two different functionals --- are discussed.


\end{abstract}

\vspace{-0.5cm}

\section{Introduction} \label{Sec: Introduction}

Kumar et al. (2018a; 2020) and  Kumar and Lopez-Pamies (2020) have recently established that any attempt at the formulation of a complete macroscopic theory of fracture nucleation and propagation in nominally elastic brittle solids must account for three intrinsic material properties: ($i$) the elasticity of the solid, ($ii$) its strength, and ($iii$) its critical energy release rate. 

In the most basic setting,  that of an isotropic linearly elastic brittle solid occupying an open bounded domain $\mathrm{\Omega}\subseteq \mathbb{R}^3$ that is subjected to monotonic and quasistatic loading conditions\footnote{We will come back to the more general case of non-monotonic and dynamic loading conditions in Section \ref{Sec: Final comments} below.} in time $t\in[0,T]$, these properties are characterized by:
\begin{enumerate}[]

\item{the stored-energy function
\begin{equation}
W(\bfE)=\mu\,{\rm tr}\,\bfE^2+\dfrac{\lambda}{2}\left({\rm tr}\,\bfE\right)^2\label{Ingredient1 W}
\end{equation}
describing elasticity of the solid, that is, the elastic energy (per unit undeformed volume) stored by the solid when deformed a strain $\bfE$,}

\item{the strength surface 
\begin{equation}
\mathcal{F}(\bfsig)=0\label{Ingredient2 Fsigma}
\end{equation}
describing the strength of the solid, that is, the set of critical stresses at which the solid fractures when subjected to monotonically increasing uniform stress $\bfsig$, and}

\item{the critical energy release rate
\begin{equation}
G_c\label{Ingredient3 Gc}
\end{equation}
describing the intrinsic fracture energy of the solid, that is, the amount of energy (per unit undeformed area) required to create new surface in the solid from an existing crack.}

\end{enumerate}

In contrast to the specific functional forms (quadratic in $\bfE$ and constant) of the stored-energy function (\ref{Ingredient1 W}) and critical energy release rate (\ref{Ingredient3 Gc}), there is \emph{no} one-size-fits-all functional form for the strength surface (\ref{Ingredient2 Fsigma}), since different solids can feature very different strength surfaces. In this work, for definiteness, we shall restrict attention to strength surfaces of the Drucker-Prager form   
\begin{equation}
\mathcal{F}(\bfsig)\equiv\sqrt{\mathcal{J}_2}+\dfrac{\sts}
{\sqrt{3}\left(3 \shs-\sts\right)} \mathcal{I}_1-\dfrac{\sqrt{3}\shs \sts}
{3\shs-\sts}=0,\label{DP}
\end{equation}
where $\mathcal{I}_1={\rm tr}\,\bfsig$ and $\mathcal{J}_2=\frac{1}{2}\,{\rm tr}\,\bfsig^2_{D}$ stand for the first and second principal invariants of the stress $\bfsig$ and the deviatoric stress $\bfsig_D=\bfsig-\frac{1}{3}({\rm tr}\,\bfsig)\bfI$ tensors, while the material constants $\sts>0$ and $\shs>0$ denote the uniaxial tensile and hydrostatic strengths of the solid, that is, they denote the critical stress values at which fracture nucleates under states of monotonically increased uniaxial tension $\bfsig={\rm diag}(\sigma>0,0,0)$ and tensile hydrostatic stress $\bfsig={\rm diag}(\sigma>0,\sigma>0,\sigma>0)$, respectively.
\begin{remark}
According to our choice of signs in (\ref{DP}), any stress state such that $$\left(3 \shs-\sts\right)\mathcal{F}(\bfsig)\geq0$$ is in violation of the strength of the solid. The strength of the vast majority of solids, especially hard solids (e.g., glass), is such that $3 \shs-\sts>0$. Nonetheless, $3 \shs-\sts<0$ for some soft solids (e.g., natural rubber).
\end{remark}
\begin{remark}
The two-material-parameter strength surface (\ref{DP}), originally introduced by Drucker and Prager (1952) to model the yielding of soils, is arguably the simplest model that has proven capable of describing reasonably well the strength of many nominally brittle solids (Lopez-Pamies, 2023; Kumar et al., 2020; 2022; 2024; Kamarei et al., 2024), thus its use here as a representative template.  
\end{remark}
\begin{remark}
For given uniaxial tensile and hydrostatic strengths $\sts$ and $\shs$, the strength surface (\ref{DP}) predicts the shear, biaxial tensile, and uniaxial compressive strengths
\begin{equation*}
\sss=\dfrac{1}{\sqrt{3}}\left(-\dfrac{1}{3}+\dfrac{\shs}{\sts}\right)^{-1}\shs,\quad \sbs=\left(\dfrac{1}{3}+\dfrac{\shs}{\sts}\right)^{-1}\shs,\quad {\rm and}\quad \scs=\left(-\dfrac{2}{3}+\dfrac{\shs}{\sts}\right)^{-1}\shs,
\end{equation*}
which are defined as the critical stress values at which fracture nucleates under uniform states of monotonically increased shear stress $\bfsig= {\rm diag}(\sigma>0, -\sigma, 0)$, biaxial tension $\bfsig={\rm diag}(\sigma>0,\sigma>0,0)$, and uniaxial compression $\bfsig={\rm diag}(\sigma<0,0,0)$. Direct use of these relations (or analogous ones for any other multi-axial stress state of interest) allows to rewrite (\ref{DP}) in terms of different pairs of critical strength constants. For hard solids, it is customary to use $\sts$ and $\scs$, while for soft solids it is more convenient to use $\sts$ and $\shs$. In this work, as indicated by (\ref{DP}), we favor the latter parametrization.
\end{remark}

\subsection{The problem}

Consider hence a body made of an isotropic linearly elastic brittle solid with stored-energy function (\ref{Ingredient1 W}), Drucker-Prager strength surface (\ref{DP}), and critical energy release rate (\ref{Ingredient3 Gc}) that, initially, at time $t=0$, occupies the open bounded domain $\mathrm{\Omega}$. We denote the boundary of the body by $\partial\mathrm{\Omega}$, its outward unit normal by $\bfN$, and identify material points in the body by their initial position vector
\begin{equation*}
\bfX\in \mathrm{\Omega}.
\end{equation*}

The body is subjected to a body force (per unit undeformed volume) $\textbf{b}(\bfX,t)$, a displacement $\overline{\bfu}(\bfX,t)$ on a part $\partial \mathrm{\Omega}_{\mathcal{D}}$ of the boundary, and a surface force (per unit undeformed area) $\overline{\bfs}(\bfX,t)$ on the complementary part $\partial\mathrm{\Omega}_{\mathcal{N}}=\partial\mathrm{\Omega}\setminus\partial \mathrm{\Omega}_{\mathcal{D}}$. In response to these stimuli --- all of which, again, are assumed to be applied monotonically and quasistatically in time --- the position vector $\bfX$ of a material point in the body will move to a new position specified by 
\begin{equation}
\bfx=\bfX+\bfu(\bfX,t),\label{def-mapping}
\end{equation}
where $\bfu(\bfX,t)$ is the displacement field. We write the associated strain at $\bfX$ and $t$ as
\begin{equation*}
\bfE(\bfu)=\dfrac{1}{2}\left(\nabla\bfu+\nabla\bfu^T\right).
\end{equation*}

In addition to the deformation (\ref{def-mapping}), the applied body force and boundary conditions may result in the nucleation and subsequent propagation of cracks in the body. We describe such cracks in a regularized fashion via the phase field
\begin{equation*}
v=v(\bfX,t)
\end{equation*}
taking values in the range $[0, 1]$. 

\subsection{The phase-field fracture theory of Kumar, Francfort and Lopez-Pamies (2018)}

According to the phase-field fracture formulation put forth by Kumar, Francfort and Lopez-Pamies (2018), the displacement field $\bfu_k(\bfX)=\bfu(\bfX,t_k)$ and phase field $v_k(\bfX)=v(\bfX,t_k)$ at any material point $\bfX \in \overline{\mathrm{\Omega}}=\mathrm{\Omega}\cup\partial\mathrm{\Omega}$ and at any given discrete time $t_k\in\{0=t_0,t_1,...,t_m,t_{m+1},...,t_M=T\}$ are determined by the system of coupled partial differential equations (PDEs)
\begin{equation}
\left\{\begin{array}{ll}
{\rm Div}\left[(v_{k}^2+\eta_{\varepsilon})\dfrac{\partial W}{\partial \bfE}(\bfE(\bfu_{k}))\right]+\textbf{b}(\bfX,t_k)=\textbf{0},& \,\bfX\in\mathrm{\Omega}\vspace{0.2cm}\\
\bfu_k(\bfX)=\overline{\bfu}(\bfX,t_k), & \, \bfX\in\partial\mathrm{\Omega}_{\mathcal{D}}\vspace{0.2cm}\\
\left[(v_{k}^2+\eta_{\varepsilon})\dfrac{\partial W}{\partial\bfE}(\bfE(\bfu_k))\right]\bfN=\overline{\textbf{s}}(\bfX,t_k),& \, \bfX\in\partial\mathrm{\Omega}_{\mathcal{N}}
\end{array}\right. \label{BVP-u-theory}
\end{equation}
and
\begin{equation}
\left\{\begin{array}{lll}
\varepsilon\, \delta^\varepsilon G_c \bigtriangleup v_k=\dfrac{8}{3}v_{k} W(\bfE(\bfu_k))+\dfrac{4}{3}c_\texttt{e}(\bfX,t_{k})-\dfrac{\delta^\varepsilon G_c}{2\varepsilon},&\, \mbox{ if } v_{k}(\bfX)< v_{k-1}(\bfX),& \bfX\in \mathrm{\Omega}\vspace{0.2cm}\\
\varepsilon\, \delta^\varepsilon G_c \bigtriangleup v_k\geq\dfrac{8}{3}v_{k} W(\bfE(\bfu_k))+\dfrac{4}{3}c_\texttt{e}(\bfX,t_{k})-\dfrac{\delta^\varepsilon G_c}{2\varepsilon},&\, \mbox{ if } v_{k}(\bfX)=1\; \mbox{ or }\; v_{k}(\bfX)= v_{k-1}(\bfX)>0, & \bfX\in \mathrm{\Omega}\vspace{0.2cm}\\
v_k(\bfX)=0,&\, \mbox{ if } v_{k-1}(\bfX)=0, & \bfX\in \mathrm{\Omega}\vspace{0.2cm}\\
\nabla v_k\cdot\bfN=0,&  &\, \bfX\in \partial\mathrm{\Omega}
\end{array}\right.  \label{BVP-v-theory}
\end{equation}
with, e.g., initial conditions $\bfu(\bfX,0)\equiv \textbf{0}$ and $v(\bfX,0)\equiv1$, where $\nabla v_k=\nabla v(\bfX,t_k)$, $\bigtriangleup v_k=\bigtriangleup v(\bfX,t_k)$,  $\eta_{\varepsilon}$ is a positive constant of $o(\varepsilon)$, and where, as elaborated below, $c_\texttt{e}(\bfX,t)$ is a driving force whose specific constitutive prescription depends on the particular form of strength surface $\mathcal{F}(\boldsymbol{\sigma})=0$, while $\delta^\varepsilon$ is a non-negative coefficient whose specific constitutive prescription depends in turn on the particular form of $c_\texttt{e}(\bfX,t)$.

\begin{remark}\label{Remark: irreversibility}
{\rm The inequalities in (\ref{BVP-v-theory}) stem from the facts that, by definition, the phase field is bounded according to $0\leq v\leq 1$ and, by constitutive assumption, fracture is an irreversible process, in other words, healing is not allowed. 
}
\end{remark}

\begin{remark}{\rm The parameter $\varepsilon$ in (\ref{BVP-v-theory}), with units of $length$, regularizes sharp cracks. Accordingly, by definition, it can be arbitrarily small. In practice, $\varepsilon$ should be selected to be smaller than  the smallest material length scale built in (\ref{BVP-u-theory})-(\ref{BVP-v-theory}),  which comes about because of the different units of the elastic stored-energy function $W(\bfE)$ ($force/length^2$), the strength function $\mathcal{F}(\bfsig)$ ($force/length^2$), and the critical energy release rate $G_c$ ($force/length$); see, e.g., Appendix C in (Kumar et al., 2024) and Appendix B in (Kamarei et al., 2024). As a rule of thumb, it typically suffices to set $\varepsilon<3G_c/(16\mathcal{W}_{\texttt{ts}})$, where $\mathcal{W}_{\texttt{ts}}$ is given by expression (\ref{Wts})$_1$ below.
}\label{Remark_length}
\end{remark}

\begin{remark}{\rm With regard to the preceding remark, it is important to emphasize that models of fracture in which $\varepsilon$ --- or any such type of a length scale parameter --- cannot be taken arbitrarily small are \emph{not} approximations of any sharp fracture model. Interestingly, by virtue of the finite value of $\varepsilon$ that may be chosen in their implementation, those models may feature an apparent strength. That strength is an artifact of such an $\varepsilon$, one that disappears as $\varepsilon\searrow 0$, and \emph{not} actual material strength. 
}\label{Fixed epsilon}
\end{remark}

\begin{remark}{\rm On their own, the governing equations (\ref{BVP-u-theory}) and (\ref{BVP-v-theory}) are standard second-order PDEs for the displacement field $\bfu_k(\bfX)$ and the phase field $v_k(\bfX)$, the latter of which is additionally subjected to standard variational inequalities. Accordingly, their numerical solution is amenable to a finite element (FE) staggered scheme in which (\ref{BVP-u-theory}) and (\ref{BVP-v-theory}) are discretized with finite elements and solved iteratively one after the other at every time step $t_k$ until convergence is reached. In the next section, as the main result of this Note, we show that the solution pair $(\bfu^{\varepsilon}_k,v^{\varepsilon}_k)$ computed in such a staggered approach corresponds in fact to the fields that minimize separately two different functionals. 
}\label{Numerics}
\end{remark}

\paragraph{The driving force $c_{\emph{\texttt{e}}}(\bfX,t)$ and coefficient $\delta^{\varepsilon}$.}  For a solid whose strength is characterized by the Drucker-Prager strength surface (\ref{DP}), Kumar et al. (2020; 2022) and Kamarei et al. (2024) have worked out different constitutive prescriptions for the driving force $c_\texttt{e}(\bfX,t)$ and coefficient $\delta^{\varepsilon}$ that are equivalent in the limit as $\varepsilon\searrow 0$, but contain different corrections of $O(\varepsilon^0)$. Here, for definiteness, we consider the constitutive prescription proposed by Kamarei et al. (2024). It reads
\begin{align}
c_{\texttt{e}}(\bfX,t)=\alpha_2\sqrt{\mathcal{J}_2}+
\alpha_1 \mathcal{I}_1\qquad {\rm and}\qquad \delta^\varepsilon=f\left(\dfrac{\shs}{\sts}\right)\dfrac{3 G_c}{16\mathcal{W}_{\texttt{ts}}\varepsilon}+g\left(\dfrac{\shs}{\sts}\right),\label{cehat}
\end{align}
where
\begin{align*}
\alpha_1=\dfrac{1}{\shs}\delta^\varepsilon\dfrac{G_c}{8\varepsilon}-\dfrac{2\mathcal{W}_{\texttt{hs}}}{3\shs},\qquad \alpha_2=\dfrac{\sqrt{3}(3\shs-\sts)}{\shs\sts}\delta^\varepsilon\dfrac{G_c}{8\varepsilon}+
\dfrac{2\mathcal{W}_{\texttt{hs}}}{\sqrt{3}\shs}-\dfrac{2\sqrt{3}\mathcal{W}_{\texttt{ts}}}{\sts},%\label{alphas}
\end{align*}
\begin{align}
\mathcal{W}_{\texttt{ts}}=\dfrac{\sts^2}{2E}, \qquad \mathcal{W}_{\texttt{hs}}=\dfrac{\shs^2}{2\kappa},\label{Wts}
\end{align}
\begin{equation*}
\left\{\begin{array}{l}
\mathcal{I}_1={\rm tr}\,\bfsig=3\kappa v^2 {\rm tr}\,\bfE\vspace{0.2cm}\\
\mathcal{J}_2=\dfrac{1}{2}\,{\rm tr}\,\bfsig^2_D=2\mu^2 v^4 {\rm tr}\,\bfE^2_D, \quad \bfE_D=\bfE-\dfrac{1}{3}\left({\rm tr}\,\bfE\right)\bfI,
\end{array}\right.,
\end{equation*}
and where we have made use of the classical connections $E=\mu(3\lambda+2\mu)/(\lambda+\mu)$ and $\kappa=\lambda+\frac{2}{3}\mu$ between the Young's modulus $E$ and bulk modulus $\kappa$ and the Lam\'e constants $\mu$ and $\lambda$; the specific forms of the non-negative functions $f$ and $g$ in (\ref{cehat})$_2$ can be found in Section 3.4 of (Kamarei et al., 2024).

The constitutive prescription (\ref{cehat}) for the driving force $c_{\texttt{e}}$ and coefficient $\delta^\varepsilon$ leads to a complete macroscopic theory of fracture (\ref{BVP-u-theory})-(\ref{BVP-v-theory}), one that over the past lustrum has been validated through direct comparisons with experiments on a wide range of nominally elastic brittle solids under a broad range of loading conditions (Kumar et al., 2018a; 2018b; 2020; 2022; 2024; Kumar and Lopez-Pamies, 2021; Kamarei et al., 2024). 

The central reason for the apparent status of the PDEs (\ref{BVP-u-theory})-(\ref{BVP-v-theory}) as a complete theory of fracture nucleation and propagation in a nominally elastic brittle solid subjected to monotonic and quasistatic loading conditions is that, by construction, they encode the competition between the elastic energy (\ref{Ingredient1 W}), the strength surface (\ref{Ingredient2 Fsigma}), and the critical energy release rate (\ref{Ingredient3 Gc}) of the solid in a way that is consistent with the large body of experimental evidence on fracture that has been amassed since the end of the 1800s.

The main objective of this Note is to show that the competition described by the PDEs (\ref{BVP-u-theory})-(\ref{BVP-v-theory}) can be recast variationally as the minimization of two different functionals. We do so in Section \ref{Sec: Variational}. The resulting variational structure happens to be of the same type that has recently allowed to prove existence of solutions in the variational approach to sharp Griffith fracture (Francfort and Marigo, 1998) accounting for boundary loads (Larsen, 2021). It is also of the same type of variational structure that describes the original phase-field theory without material strength when implemented as an alternating minimization (Bourdin et al., 2000; 2006). We devote Section \ref{Sec: Original Phase-Field} to explaining this overarching connection. There, we also point to some misconceptions in the literature on the interpretation of the original phase-field theory as typically implemented. We conclude in Section \ref{Sec: Final comments} by recording a number of final comments.

\section{The variational formulation of the phase-field fracture theory (\ref{BVP-u-theory})-(\ref{BVP-v-theory})} \label{Sec: Variational}

\subsection{The deformation energy functional $\mathcal{E}^{\varepsilon}_{d}(\bfu_k;v_k)$} 

For any fixed phase field $v_k(\bfX)\in[0,1]$, a standard calculation suffices to show that equations (\ref{BVP-u-theory}) are nothing more than the Euler-Lagrange equations associated with variations in 
\begin{equation*}
\bfu_k\in\mathcal{A}_u=\left\{\bfu_k\in H^{1}(\mathrm{\Omega};\mathbb{R}^3): \bfu_k(\bfX)=\overline{\bfu}(\bfX,t_k),\;\bfX\in\partial\mathrm{\Omega}_{\mathcal{D}}\right\}
\end{equation*}
of the \emph{deformation energy functional}
\begin{equation}
\mathcal{E}^{\varepsilon}_{d}(\bfu_k;v_k):=\displaystyle\int_{\mathrm{\Omega}} (v_k^2+\eta_{\varepsilon}) W(\bfE(\bfu_k)){\rm d}\bfX-\displaystyle\int_{\mathrm{\Omega}}\textbf{b}(\bfX,t_k)\cdot\bfu_k\,{\rm d}\bfX-\displaystyle\int_{\partial\mathrm{\Omega}_{\mathcal{N}}}\overline{\textbf{s}}(\bfX,t_k)\cdot\bfu_k\,{\rm d}\bfX.\label{Ed}
\end{equation}
Assuming that $\textbf{b}(\bfX,t_k)\in L^2(\mathrm{\Omega};\mathbb{R}^3)$, $\overline{\textbf{\bfu}}(\bfX,t_k)\in L^{2}(\partial\mathrm{\Omega}_{\mathcal{D}};\mathbb{R}^3)$, and $\overline{\textbf{s}}(\bfX,t_k)\in L^{2}(\partial\mathrm{\Omega}_{\mathcal{N}};\mathbb{R}^3)$, and noting that $\eta_{\varepsilon}>0$, we can invoke the classical theorems of existence and uniqueness in linear elastostatics (see, e.g., Fichera, 1973) to readily establish that 
\begin{equation}
(\bfu^\varepsilon_k;v_k)=\underset{\begin{subarray}{c}
  \bfu_k\in\mathcal{A}_u
  \end{subarray}}{\arg\min}\, \mathcal{E}^{\varepsilon}_{d}(\bfu_k;v_k)\label{Minimizer-u}
\end{equation}
exists, is unique, and is the solution of equations (\ref{BVP-u-theory}).

\subsection{The fracture functional $\mathcal{E}^{\varepsilon}_{f}(v_k;\bfu_k)$} 

Define the ``undamaged'' driving force 
\begin{equation*}
\widehat{c}_{\texttt{e}}(\bfE(\bfu))=\mu\alpha_2\sqrt{2\,{\rm tr}\,\bfE^2_D(\bfu)}+
3\kappa \alpha_1  {\rm tr}\,\bfE(\bfu)\, \textrm{ so that }\, c_{\texttt{e}}(\bfX,t)=v^2\widehat{c}_{\texttt{e}}(\bfE(\bfu)).%\label{chat}
\end{equation*}
In terms of this driving force, equations (\ref{BVP-v-theory}) can be rewritten as
\begin{equation}
\left\{\begin{array}{lll}
\varepsilon\, \delta^\varepsilon G_c \bigtriangleup v_k=\dfrac{8}{3}v_{k} W(\bfE(\bfu_k))+\dfrac{4}{3}v_k^2\widehat{c}_\texttt{e}(\bfE(\bfu_k))-\dfrac{\delta^\varepsilon G_c}{2\varepsilon},&\, \mbox{ if } v_{k}(\bfX)< v_{k-1}(\bfX),& \bfX\in \mathrm{\Omega}\vspace{0.2cm}\\
\varepsilon\, \delta^\varepsilon G_c \bigtriangleup v_k\geq\dfrac{8}{3}v_{k} W(\bfE(\bfu_k))+\dfrac{4}{3}v_k^2\widehat{c}_\texttt{e}(\bfE(\bfu_k))-\dfrac{\delta^\varepsilon G_c}{2\varepsilon},&\, \mbox{ if } v_{k}(\bfX)=1\; \mbox{ or }\; v_{k}(\bfX)= v_{k-1}(\bfX)>0, & \bfX\in \mathrm{\Omega}\vspace{0.2cm}\\
v_k(\bfX)=0,&\, \mbox{ if } v_{k-1}(\bfX)=0, & \bfX\in \mathrm{\Omega}\vspace{0.2cm}\\
\nabla v_k\cdot\bfN=0,&  &\, \bfX\in \partial\mathrm{\Omega} .
\end{array}\right.  \label{BVP-v-theory-chat}
\end{equation}
For any fixed regularization length $\varepsilon>0$ and any fixed displacement field $\bfu_k(\bfX)\in H^{1}(\mathrm{\Omega};\mathbb{R}^3)$, yet another standard calculation suffices
to show that equations (\ref{BVP-v-theory-chat}) are the Euler-Lagrange equations associated with variations in
\begin{equation}
v_k\in\mathcal{A}_v=\left\{v_k\in H^{1}(\mathrm{\Omega};\mathbb{R}): 0\leq v_k\leq 1,\, v_k\leq v_{k-1} \right\}\label{Av}
\end{equation}
of the \emph{fracture functional}
\begin{equation}
\mathcal{E}^{\varepsilon}_{f}(v_k;\bfu_k):=\displaystyle\int_{\mathrm{\Omega}} v_k^2 W(\bfE(\bfu_k)){\rm d}\bfX+\displaystyle\int_{\mathrm{\Omega}}\dfrac{v_k^3}{3}\,\widehat{c}_{\texttt{e}}(\bfE(\bfu_k)){\rm d}\bfX+\dfrac{3\,\delta^\varepsilon G_c}{8} \int_{\mathrm{\Omega}}\left(\dfrac{1-v_k}{\varepsilon}+\varepsilon\nabla v_k\cdot\nabla v_k\right)\,{\rm d}\bfX.\label{Ef}
\end{equation}
Since, as outlined next, $\mathcal{E}^{\varepsilon}_{f}$ has minimizers in the admissible set (\ref{Av}), we have that 
\begin{equation}
(v^\varepsilon_k;\bfu_k)\in\underset{\begin{subarray}{c}
  v_k\in\mathcal{A}_v
  \end{subarray}}{\arg\min}\, \mathcal{E}^{\varepsilon}_{f}(v_k;\bfu_k)\label{Minimizer-v}
\end{equation}
exists and is a solution of equations (\ref{BVP-v-theory-chat}).

We now briefly describe the existence of minimizers of $\mathcal{E}^{\varepsilon}_{f}$. We suppose that $v_{k-1}=1$ for simplicity, the general case is similar. The idea is a straightforward application of the Direct Method in the Calculus of Variations. We consider a sequence $v_n$ such that $\mathcal{E}^{\varepsilon}_{f}(v_n;\bfu_k)\rightarrow \inf \mathcal{E}^{\varepsilon}_{f}(\cdot;\bfu_k)$. Note that
\[
-\frac13\displaystyle\int_{\mathrm{\Omega}}\,|\widehat{c}_{\texttt{e}}(\bfE(\bfu_k))|{\rm d}\bfX+
\dfrac{3\,\delta^\varepsilon G_c}{8} \int_{\mathrm{\Omega}}\varepsilon\nabla v_n\cdot\nabla v_n\,{\rm d}\bfX
\leq
\mathcal{E}^{\varepsilon}_{f}(v_n;\bfu_k)\leq
\displaystyle\int_{\mathrm{\Omega}}  (W(\bfE(\bfu_k))+\,|\widehat{c}_{\texttt{e}}(\bfE(\bfu_k))|){\rm d}\bfX<\infty.
\]
It follows that $\{\nabla v_n\}$ is bounded in $L^2(\mathrm{\Omega};\mathbb{R}^3)$, and since $|v_n|\leq 1$, $\{v_n\}$ is bounded in $H^1(\mathrm{\Omega};\mathbb{R})$. Therefore, a subsequence, not relabeled, satisfies $v_n\wc v_\infty$ in $H^1(\mathrm{\Omega};\mathbb{R})$ for some $v_\infty\in H^1(\mathrm{\Omega};\mathbb{R})$, and so $v_n\rightarrow v_\infty$ strongly in $L^2(\mathrm{\Omega};\mathbb{R})$.  As $\nabla v_n\wc \nabla v_\infty$ in $L^2(\mathrm{\Omega};\mathbb{R}^3)$, we get
\[
\int_{\mathrm{\Omega}}\left(\varepsilon\nabla v_\infty\cdot\nabla v_\infty\right)\,{\rm d}\bfX\leq \liminf_{n\rightarrow \infty}
\int_{\mathrm{\Omega}}\left(\varepsilon\nabla v_n\cdot\nabla v_n\right)\,{\rm d}\bfX
\]
by the convexity of this term.

The first two terms in $\mathcal{E}^{\varepsilon}_{f}$ are of the form $\int_{\mathrm{\Omega}}v_k^p F\,{\rm d}\bfX$, where $F$ is integrable and $p=2,3$.  Since $v_n\rightarrow v_\infty$ strongly in $L^2(\mathrm{\Omega};\mathbb{R})$, for a subsequence, not relabeled, $v_n\rightarrow v_\infty$ a.e. Therefore, $v_n^p F\rightarrow v_\infty^p F$ a.e., and since $v_n\in [0,1]$, we have $v_n^p|F|\leq |F|$. The Dominated Convergence Theorem then implies
\[
\int_{\mathrm{\Omega}}v_n^pF\,{\rm d}\bfX\rightarrow \int_{\mathrm{\Omega}}v_\infty^pF\,{\rm d}\bfX.
\]
The term in $(1-v_k)$ converges similarly.

Putting all the terms together, we get
\[
\mathcal{E}^{\varepsilon}_{f}(v_\infty;\bfu_k)\leq \liminf_{n\rightarrow \infty}\mathcal{E}^{\varepsilon}_{f}(v_n;\bfu_k)= \inf \mathcal{E}^{\varepsilon}_{f}(\cdot;\bfu_k),
\]
and so $v_\infty$ is a minimizer. Finally, note that for any irreversibility constraint (of the type included in (\ref{Av})) such that  $v_n\leq V$ a.e. for some $V$, then since $v_n\rightarrow v_\infty$ a.e., we also get $v_\infty\leq V$ a.e.

\begin{remark}
Contrary to the uniqueness of the minimizer (\ref{Minimizer-u}) of the deformation energy functional (\ref{Ed}), the minimizer (\ref{Minimizer-v}) of the fracture functional (\ref{Ef}) need \emph{not} be unique.
\end{remark}

\subsection{The variational principle} 

Having defined the deformation energy and fracture functionals (\ref{Ed}) and (\ref{Ef}) and having established that they have minimizers, we are now ready to recast the phase-field fracture theory described by the PDEs (\ref{BVP-u-theory})-(\ref{BVP-v-theory}) as a variational principle. 

The idea --- as in the standard implementation of the original phase-field fracture theory (see below) --- is to alternately minimize the deformation energy and fracture functionals (\ref{Ed}) and (\ref{Ef}) at every time step $t_k$ until reaching a converged solution pair $(\mathfrak{u}_k^{\varepsilon},\mathfrak{v}_k^{\varepsilon})$. Doing so leads to sequences $\bfu_i$, $v_i$ satisfying
\[
\bfu_i\mbox{ minimizes } \mathcal{E}^{\varepsilon}_{d}(\cdot\;; v_{i-1})
\]
and
\[
v_i\mbox{ minimizes } \mathcal{E}^{\varepsilon}_{f}(\cdot\;; \bfu_{i}).
\]
The expectation --- as, again, in the standard implementation of the original phase-field fracture theory (see below) --- is that the limits $\mathfrak{u}_k^{\varepsilon}$ and $\mathfrak{v}_k^{\varepsilon}$ of $\bfu_i$ and $v_i$, respectively, satisfy the minimality conditions
\[
\mathfrak{u}_k^{\varepsilon}\mbox{ minimizes } \mathcal{E}^{\varepsilon}_{d}(\cdot\;; \mathfrak{v}_k^{\varepsilon})
\]
and
\[
\mathfrak{v}_k^{\varepsilon}\mbox{ minimizes } \mathcal{E}^{\varepsilon}_{f}(\cdot\;; \mathfrak{u}_k^{\varepsilon})
\]
and therefore satisfy the PDEs (\ref{BVP-u-theory})-(\ref{BVP-v-theory}). At present, there is no mathematical proof. However, as noted in the Introduction, the verity of this expectation has been supported by numerous simulations of experiments on a wide range of nominally elastic brittle solids under a broad range of loading conditions, which have also served as validation results for the theory. In the sequel, we provide some insight into why this is the case.

\begin{remark}
Larsen (2021) has recently shown that a similar variational principle based on the minimization of two different functionals in the setting of sharp --- as opposed to phase-field --- fracture is necessary in order to have existence of solutions in the variational approach to fracture accounting for boundary loads. 
\end{remark}

\subsection{Some remarks on the competitions set up by the variational principle} 

On the one hand, the minimization of the deformation energy functional (\ref{Ed}) is nothing more than the classical variational statement that forces in the body satisfy balance of linear and angular momenta for the specific case when the body is made of a linear elastic solid.

On the other hand, the minimization of the fracture functional (\ref{Ef}) states that whether cracks nucleate and/or propagate depends on a competition among all three intrinsic material properties (\ref{Ingredient1 W})-(\ref{Ingredient3 Gc}) of the solid at states when the body is in elastostatic equilibrium. 

\paragraph{Uniform fields.} In particular, for the basic case when the body is subjected to a \emph{uniform} strain $\bfE(\bfu_k)=\overline{\bfE}_k$, with help of the notation 
\begin{align*}
\overline{\bfS}_k=\dfrac{\partial W}{\partial\bfE}(\overline{\bfE}_k)
\end{align*}
for the corresponding uniform ``undamaged'' stress, noting that 
\begin{align*}
\widehat{c}_{\texttt{e}}(\overline{\bfE}_k)=\dfrac{3 \delta^\varepsilon G_c}{8\varepsilon}\dfrac{(3\shs-\sts)}{\sqrt{3}\shs\sts}\left(\mathcal{F}(\overline{\bfS}_k)+\dfrac{\sqrt{3}\shs\sts}{3\shs-\sts}\right) +O(\varepsilon^0),
\end{align*}
where we recall that the strength function $\mathcal{F}$ is defined in (\ref{DP}), the fracture functional (\ref{Ef}) specializes to 
\begin{align}
\mathcal{E}^{\varepsilon}_{f}(v_k;\bfu_k)=&\displaystyle\int_{\mathrm{\Omega}} v_k^2 W(\overline{\bfE}_k){\rm d}\bfX+\displaystyle\int_{\mathrm{\Omega}}\dfrac{v_k^3}{3}\left[\dfrac{3 \delta^\varepsilon G_c}{8\varepsilon}\dfrac{(3\shs-\sts)}{\sqrt{3}\shs\sts}\left(\mathcal{F}(\overline{\bfS}_k)+\dfrac{\sqrt{3}\shs\sts}{3\shs-\sts}\right) +O(\varepsilon^0)\right]{\rm d}\bfX+\vspace{0.2cm}\nonumber\\ 
&\dfrac{3\,\delta^\varepsilon G_c}{8} \int_{\mathrm{\Omega}}\left(\dfrac{1-v_k}{\varepsilon}+\varepsilon\nabla v_k\cdot\nabla v_k\right)\,{\rm d}\bfX.\label{Ef-uniform}
\end{align}
For sufficiently small regularization length $\varepsilon$, in regions where $v_k=1$, the only terms that compete in the minimization of this functional are the second (i.e., the strength) and the third (i.e., the fracture energy) integrals. The first integral (i.e., the elastic energy) is inconsequential in this case. 

In view of the definition of the strength function $\mathcal{F}$ in (\ref{DP}), the second integral in (\ref{Ef-uniform}) starts at a value of 0 when $\overline{\bfE}_k=\textbf{0}$. As the strain $\overline{\bfE}_k$ deviates from $\textbf{0}$ along the given loading path, so does the second integral, which may become large enough to compete with the third integral and the minimizing displacement (\ref{Minimizer-u}) to force the localization of the phase field $v_k$ near $v_k=0$ and hence the nucleation of fracture. As this localization of $v_k$ occurs, the first integral in (\ref{Ef-uniform}) may become of comparable order to the other two integrals in a way that further enhances localization. Consistent with experimental observations, the numerical experiments referenced in the Introduction indicate that this localization of $v_k$ (when the strain is initially uniform in the body) happens when 
\begin{align*}
\mathcal{F}(\overline{\bfS}_k)=0,
\end{align*}
that is, when the strength surface (\ref{DP}) of the solid is first violated. 

Importantly --- in stark contrast to the first integral which is non-negative --- the second integral in (\ref{Ef-uniform}) can be positive or negative. When it is negative, or when it is positive but not sufficiently large, there is no incentive to localize the phase field $v_k$ and hence the nucleation of fracture does not occur. Consistent with experimental observations, the numerical experiments referenced in the Introduction indicate that fracture nucleation does \emph{not} occur so long as 
\begin{align*}
(3\shs-\sts)\mathcal{F}(\overline{\bfS}_k)\leq 0.
\end{align*}

\paragraph{Non-uniform fields due to the presence of large cracks.} For the opposite basic case when the body contains large\footnote{``Large'' refers to large relative to the characteristic size of the underlying heterogeneities in the solid under investigation. By the same token, ``small'' refers to sizes that are of the same order or just moderately larger than the sizes of the heterogeneities. In practice, as a rule of thumb, ``large'' refers to large relative to the material length scale $G_c/\mathcal{W}_{\texttt{ts}}$.} cracks, the numerical experiments indicate that the competition set up by the fracture functional (\ref{Ef}) describes that such cracks grow according to Griffith's criticality condition (Griffith, 1921), consistent with experimental observations.

Importantly, all three integrals in (\ref{Ef}) enter the competition that describes the growth of large cracks. In particular, because of its different sign for tensile and compressive stresses, the second integral (i.e., the strength) prevents the growth of large cracks under compressive loads, consistent with experimental observations (Liu and Kumar, 2024).

\paragraph{Arbitrary non-uniform fields.} Irrespective of the presence of large cracks, the stress field in a body is typically highly non-uniform. In this general case, consistent yet again with experimental observations, numerical experiments indicate that the competition set up by the fracture functional (\ref{Ef}) describes that fracture nucleation is governed neither solely by the strength of the solid nor solely by Griffith's criticality condition, but by an ``interpolation'' between the two. 

What is more, the results indicate that a necessary condition for the nucleation of fracture at a material point $\bfX$ is that the strength surface of the solid be violated at that point. Such a condition is \emph{not} sufficient, however, as indicated by the observation that stress singularities will always give rise to a stress state that violates the strength surface without necessarily resulting in the nucleation of fracture.

\section{A comment on the variational structure of the original formulation of Griffith phase-field fracture without material strength} \label{Sec: Original Phase-Field}

In the original phase-field fracture formulation (see, e.g., Bourdin et al., 2000; 2006), one seeks a minimizing pair $(\bfu_k^\varepsilon,v_k^\varepsilon)$ of the energy functional
\begin{equation*}
\mathcal{E}^\varepsilon(\bfu_k,v_k):=\displaystyle\int_{\mathrm{\Omega}} v_k^2 W(\bfE(\bfu_k)){\rm d}\bfX-\displaystyle\int_{\mathrm{\Omega}}\textbf{b}(\bfX,t_k)\cdot\bfu_k\,{\rm d}\bfX-\displaystyle\int_{\partial\mathrm{\Omega}_{\mathcal{N}}}\overline{\textbf{s}}(\bfX,t_k)\cdot\bfu_k\,{\rm d}\bfX + \dfrac{3\,G_c}{8} \int_{\mathrm{\Omega}}\left(\dfrac{1-v_k}{\varepsilon}+\varepsilon\nabla v_k\cdot\nabla v_k\right)\,{\rm d}\bfX.
\end{equation*}
This energy functional is \emph{not} convex and hence difficult to minimize. It is separately convex, however, and so standard implementations involve alternating minimization. Precisely, sequences of pairs $(\bfu_i, v_i)$ are found, satisfying
\[
\bfu_i\mbox{ minimizes } \mathcal{E}^{\varepsilon}(\cdot\;, v_{i-1})
\]
and
\[
v_i\mbox{ minimizes } \mathcal{E}^{\varepsilon}(\bfu_{i}, \;\cdot).
\]
Looking at the terms that are independent of $v_k$ and $\bfu_k$ respectively, this is equivalent to saying
\[
\bfu_i\mbox{ minimizes } \mathcal{E}^{\varepsilon}_{d}(\cdot\;; v_{i-1})
\]
and
\[
v_i\mbox{ minimizes } \mathcal{E}^{\varepsilon}_{G}(\cdot\;; \bfu_{i}),
\]
where we recall that the deformation energy functional $\mathcal{E}^{\varepsilon}_{d}$ is defined by (\ref{Ed}), while
\begin{equation}
\mathcal{E}^{\varepsilon}_{G}(v_k;\bfu_k):=\displaystyle\int_{\mathrm{\Omega}} v_k^2 W(\bfE(\bfu_k)){\rm d}\bfX + \dfrac{3\,G_c}{8} \int_{\mathrm{\Omega}}\left(\dfrac{1-v_k}{\varepsilon}+\varepsilon\nabla v_k\cdot\nabla v_k\right)\,{\rm d}\bfX\label{EG}
\end{equation}
defines the \emph{Griffith energy functional}; note that, as opposed to the fracture functional (\ref{Ef}), the Griffith energy functional (\ref{EG}) does \emph{not} account for the strength of the solid. As in the previous section, the expectation, without proof, is that the limits $\mathfrak{u}_k^{\varepsilon}$ and $\mathfrak{v}_k^{\varepsilon}$ of $\bfu_i$ and $v_i$, respectively, satisfy the minimality conditions
\[
\mathfrak{u}^{\varepsilon}_k\mbox{ minimizes } \mathcal{E}^{\varepsilon}_{d}(\cdot\;; \mathfrak{v}_k^{\varepsilon})
\]
and
\[
\mathfrak{v}_k^{\varepsilon}\mbox{ minimizes } \mathcal{E}^{\varepsilon}_{G}(\cdot\;; \mathfrak{u}_k^{\varepsilon}).
\]

At this stage, it is plain that the solution pairs $(\bfu_k^{\varepsilon},v_k^{\varepsilon})$ described by the phase-field fracture theory of Kumar, Francfort, and Lopez-Pamies (2018) and the original phase-field fracture theory are described by the same type of variational principle. The sole difference between them is that the original phase-field fracture theory is based on the Griffith energy functional (\ref{EG}), which does \emph{not} account for the strength surface $\mathcal{F}(\bfsig)=0$ of the solid, while the phase-field fracture theory of Kumar, Francfort, and Lopez-Pamies (2018) is based on the fracture functional (\ref{Ef}), which accounts for all three material properties (\ref{Ingredient1 W})-(\ref{Ingredient3 Gc}) required for a complete theory of fracture.

\begin{remark}
We close this section by emphasizing that for neither setting, the original phase-field fracture theory and that of Kumar, Francfort, and Lopez-Pamies (2018), is there a mathematical proof (yet) that guarantees that large cracks grow according to Griffith's criticality condition. Whether they do can only be verified (thus far) numerically. For the case of the original phase-field fracture theory, this is because alternating minimization does \emph{not} necessarily lead to the global minimization of $\mathcal{E}^{\varepsilon}$. What is more, as discussed by Larsen (2023), it is possible to alter the energy functional $\mathcal{E}^{\varepsilon}$ so that both fracture nucleation and propagation are completely prevented in the alternating minimization described above, while still having $\mathcal{E}^{\varepsilon}$ $\mathrm{\Gamma}$-converge to the energy functional of sharp Griffith fracture.
\end{remark}

\section{Final comments} \label{Sec: Final comments}

In this Note, we have shown that the phase-field fracture theory of Kumar, Francfort, and Lopez-Pamies (2018) can be formulated as the minimization of two different functionals: ($i$) a deformation energy functional $\mathcal{E}_d^{\varepsilon}$ and ($ii$) a fracture functional $\mathcal{E}_f^{\varepsilon}$. While the minimization of the deformation energy functional $\mathcal{E}_d^{\varepsilon}$ determines the deformation of the body, the minimization of the fracture functional $\mathcal{E}_f^{\varepsilon}$ determines where and when cracks nucleate and propagate.

We have also shown that the original phase-field fracture theory (Bourdin et al., 2000; 2008), as typically implemented in terms of an alternating minimization procedure, is described by the same variational formulation, with the key difference that the fate of cracks is determined not by the fracture functional $\mathcal{E}_f^{\varepsilon}$, but by a Griffith energy functional $\mathcal{E}_G^{\varepsilon}$, which does \emph{not} account for the strength of solids. 

While the focus of this Note has been on the basic case of fracture in nominally linear elastic brittle solids under monotonic and quasistatic loading conditions, it is apparent that the proposed variational approach --- in terms of the minimization of two different functionals --- has the flexibility to accommodate additional physical phenomena, such as non-monotonic and dynamic loading conditions, as well as possibly the inelasticity of actual solids. From a mathematical point of view, the proposed variational approach may also provide a fruitful path to advance the analysis of fracture theories at large.

\section*{Acknowledgements}

\noindent Support for this work by the National Science Foundation through the Grants DMS--2206114 and DMS--2308169 is gratefully acknowledged. This work began during the 2024 workshop ``Fracture as an emergent phenomenon" at the Mathematisches Forschungsinstitut Oberwolfach, who we thank for hosting us.

\end{document}